\documentclass[14pt]{article}
\usepackage[english]{babel}
\usepackage{amsfonts,amssymb,amsmath,amstext,amsbsy,amsopn,amscd,amsthm,graphicx,euscript,graphicx}
\usepackage{graphics}
\newtheorem{proposition}{Proposition}
\newtheorem{corollary}{Corollary}

\theoremstyle{definition}
\newtheorem{example}{Example}
\newtheorem{remark}{Remark}

\newcommand{\eps}{\varepsilon}

\newcounter{tdfn}
\newcounter{trk}
{\endtrivlist}

\def\:{\colon}

\def\R{{\mathbb R}}
\def\Z{{\mathbb Z}}
\def\0{{\mathbf 0}}
\def\1{{\mathbf 1}}

\def\R{{\mathbb R}}

 \author{V.O.Manturov \footnote{Moscow Institute of Physics and Technology}, I.M.Nikonov\footnote{Moscow State University}}

\title{Maps from braids to virtual braids and braid representations}

\begin{document}

\maketitle

\section{Introduction}
Virtual knot theory has experienced a lot of nice features that did not appear in classical
knot theory, e.g., parity and picture-valued invariants~\cite{Man,FKMN}.

In the present paper we use virtual knot theory effects to construct new representations of classical
(pure) braids.

The paper consists of two parts: in Section 2, we use the ``$G_{n}^{3}$-approach'' by looking at triple of
collinear points, and in Section 3, we use the ``$G_{n}^{4}$-approach'' by looking at quadruples of points
belonging to the same circle or line.

\section{A new invariant of pure braids: triples of collinear points}

\subsection{The map $p_k$}

Let $\beta$ be a pure braid on $n$ strands. Then $\beta=\{\beta_k\}_{k=1}^n$ where $\beta_k\colon[0,1]\to \R^2=\mathbb C$ are the strands of the braid.

Fix an index $k\in \{1,\cdots,n\}$. The maps $p_k(\beta)_l\colon [0,1]\to S^1=\{z\in\mathbb C\mid |z|=1\}$, $l\ne k$, given by the formulas
\[
p_k(\beta)_l(t)=\frac{\beta_l(t)-\beta_k(t)}{|\beta_l(t)-\beta_k(t)|}
\]
define a braid diagram in the cylinder (we can indicate the over-under crossing structure for this projection). Denote the corresponding braid in the cylinder by $p_k(\beta)$.

\begin{proposition}
The map $p_k$ is a homomorphism from the (pure) braids on $n$ strands $PB_n$ to the group of (pure) braids in the cylinder on $n-1$ strands $CPB_{n-1}$.
\end{proposition}

Indeed, any isotopy of the braid $\beta$ induces an isotopy of the braid in the cylinder $p_k(\beta)$.

\begin{figure}[h]
\centering
\includegraphics[width=0.3\textwidth]{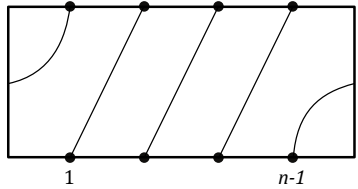}
\caption{The braid in the cylinder $\zeta$}\label{fig:zeta}
\end{figure}

In order to give an explicit description of the homomorphism $p_k$, assume that the initial points $\beta_l(0)$, $l=1,\dots,n$, lie in a circle. The Artin generators $\sigma_i$, $i=1,\dots,n-1$, of the braid group correspond to some dynamics of the points $\beta_l(0)$. Taking the compositions of these dynamics with the projection $p_k$, one gets the expressions for the corresponding braids in the cylinder (see Fig.~\ref{fig:map_pk}):

\begin{figure}[h]
\centering
\includegraphics[width=0.7\textwidth]{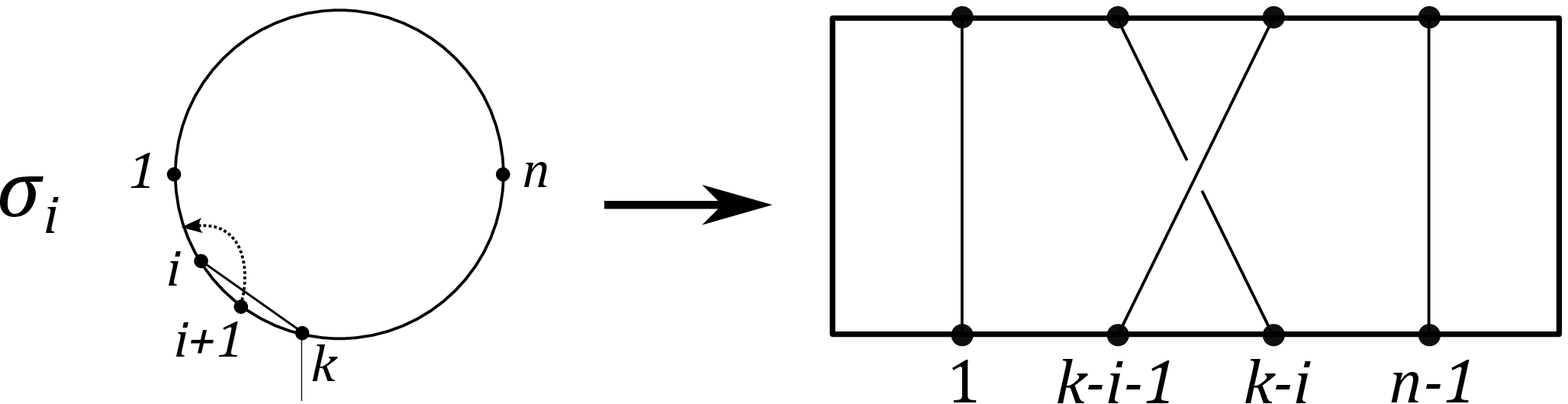}\\
\smallskip
\includegraphics[width=0.7\textwidth]{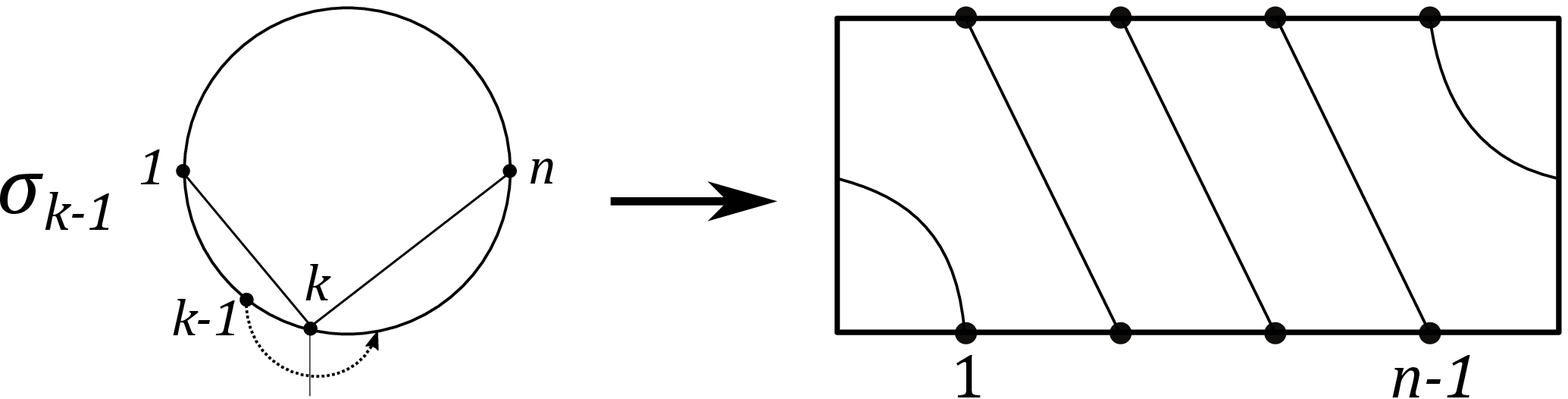}\\
\smallskip
\includegraphics[width=0.7\textwidth]{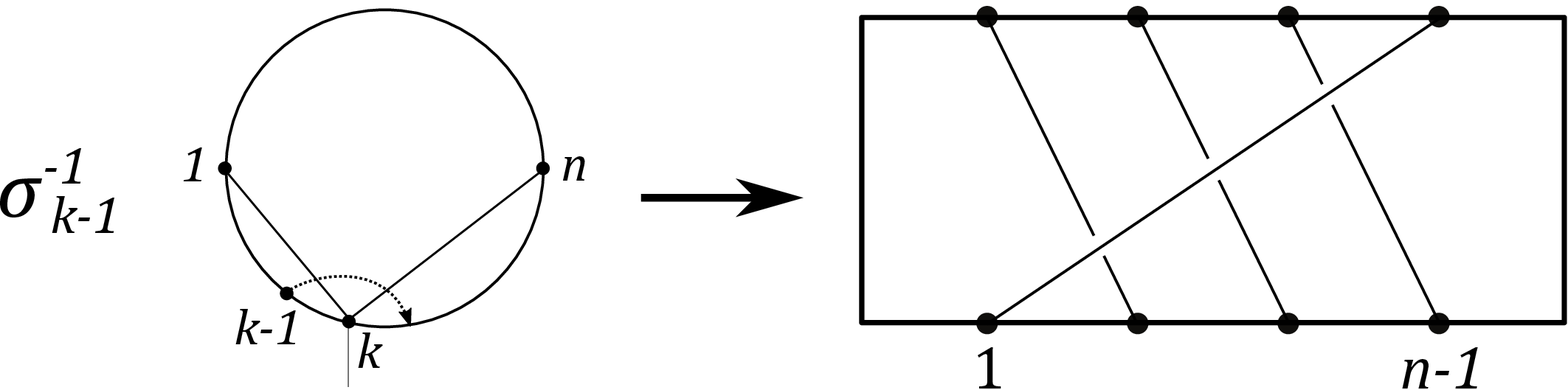}
\caption{Braids in the cylinder corresponding to the dynamics of Artin generators. The cylinder is cut into a rectangle. The cut corresponds to the downward ray from the point $k$ in the left figure.}\label{fig:map_pk}
\end{figure}

\begin{equation}\label{eq:pk}
p_k(\sigma_i^\epsilon)=\left\{
\begin{array}{cl}
\sigma_{k-i-1}^\epsilon, & i\ne k-1,k,\\
\zeta^{-1}, & i=k-1, \epsilon=1,\\
\Delta_c, & i=k-1, \epsilon=-1,\\
\Delta_c^{-1}, & i=k, \epsilon=1,\\
\zeta, & i=k, \epsilon=-1,
\end{array}
\right.
\end{equation}

where $\Delta_c=\sigma_1\cdot\sigma_2\cdots\sigma_{n-2}$ and $\zeta$ is the braid in the cylinder shown in Fig.~\ref{fig:zeta}

\subsection{The map $f_d$}

Let $\beta'$ be a pure braid diagram in the cylinder on $n'=n-1$ strands. Assume that the initial configuration of the braid $\{\beta'_k(0)\}\subset S^1$ is in general position (here this means that $\frac{\beta'_k(0)}{\beta'_l(0)}$, $k\ne l$, are transcendental numbers). Let $d\in\mathbb N$. Consider the map $f_d(z)=z^d$. Then the image $f_d(\beta')$ in $S^1$ will be a set of $n'$ strands.

Consider an intersection point (crossing) of the new braid diagram. If this crossing is the image of a crossing in $\beta'$ assign an over-under crossing structure to it so that the signs of the image crossing and the preimage crossing coincide.
If a crossing is not an image of a crossing in $\beta'$, mark it as a virtual crossing. Thus, we get a braid diagram in the cylinder with classical and virtual crossings, which we denote by $f_d(\beta')$. Equivalence classes of such diagrams modulo classical and virtual Reidemeister moves form the group of virtual braids in the cylinder $VCB_{n'}$.

\begin{proposition}
The map $f_d$ is a homomorphism from the braids in the cylinder $CB_{n'}$ to the virtual braids in the cylinder $VCB_{n'}$.
\end{proposition}

\begin{figure}[h]
\centering
\includegraphics[width=0.8\textwidth]{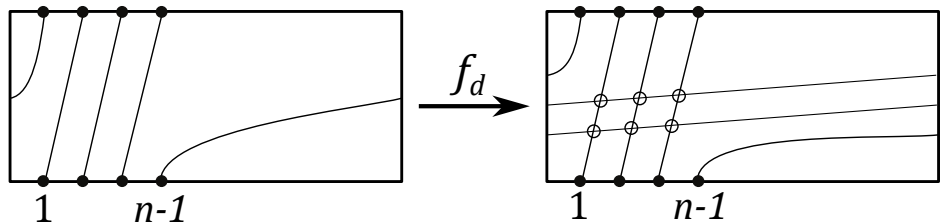}
\caption{The image of the braid $\zeta$ under the map $f_d$}\label{fig:fd_zeta}
\end{figure}

In order to give an explicit description of the homomorphism $f_d$, assume that the initial points $\beta'_l(0)$, $l=1,\dots,n-1$, lie near the point $1$ in the circle $S^1=\{z\in\mathbb C\mid |z|=1$. Then $f_d(\sigma_i)=\sigma_i$, $i=1,\dots,n-2$, because the map $f_d$ just stretches the corresponding braid diagram. On the other hand, we have (Fig.~\ref{fig:fd_zeta})
\[f_d(\zeta)=\zeta(\Delta_v\zeta)^{d-1}
\]
where $\Delta_v=\tau_1\tau_2\dots\tau_{n-2}$.

Let us define a representation of virtual cylinder braid. The group $VCB_n$ has a presentation with generators $\sigma_k$ (classical crossings) and $\tau_k$ (virtual crossings), $k=1,\dots,n$. The relations are classical and virtual Artin relations.

For a generator $\sigma_k$, we consider the matrix $\rho(\sigma_k)\in M(n)$ which coincides with the identity matrix beyond the rows and columns with numbers $k$ and $k+1$ (considered modulo $n$). In those rows and columns the matrix is equal to
\[
\left(
\begin{array}{cc}
   1-t  & t \\
    1 & 0
\end{array}
\right).
\]
Analogously, to a generator $\tau_k$ we assign the matrix $\rho(\tau_k)$ whose restriction to the rows and columns with numbers $k$ and $k+1$ is equal to
\[
\left(
\begin{array}{cc}
   0  & s \\ s^{-1} & 0
\end{array}
\right).
\]

The braid $\zeta$ corresponds to the permutation matrix
\[
\left(
\begin{array}{ccccc}
   0  & 1 & 0 & \dots &0 \\
   0  & 0 & 1 & \dots & 0 \\
   \dots  & \dots & \dots & \dots & \dots \\
   0  & 0 & \dots & 0 & 1 \\
   1  & 0 & \dots & 0 & 0
\end{array}
\right).
\]

The following proposition follows from the direct checking the relations in the virtual braid group (cf~\cite{bard}).

\begin{proposition}\label{prop:pHomomorphism}
The map $\rho$ defines a homomorphism $\rho\colon VCB_n\to GL(n,\Z[t,t^{-1},s,s^{-1}])$.
\end{proposition}

\begin{corollary}
The composition $\rho\circ f_d\circ p_k\colon PB_n\to GL(n-1,\Z[t,t^{-1},s,s^{-1}])$ is a representation of the pure braid group.
\end{corollary}

\begin{example}
Consider the element $\beta=[\psi_1{_1}\sigma_4\psi_1,\psi_2{_1}\sigma_4\sigma_3\sigma_2\sigma_1^2\sigma_2\sigma_3\sigma_4\psi_2]\in B_5$ where
$\psi_1=\sigma_3^{-1}\sigma_2\sigma_1^2\sigma_2\sigma_4^3\sigma_3\sigma_2$ and $\psi_2=\sigma_4^{-1}\sigma_3\sigma_2\sigma_1^{-2}\sigma_2\sigma_1^2\sigma_2^2\sigma_1\sigma_4^5$ which lies in the kernel of Burau representation~\cite{Big}.
Then for $t=-1, s=1$ one has
\[
\rho\circ f_2\circ p_1(\beta)=
\left(
\begin{array}{cccc}
   481  & -880 & 800 & -400  \\
   480 & -879 & 800 & -400 \\
   480 & -880 & 801 & -400\\
   480 & -880 & 800 & -399
\end{array}
\right).
\]
Thus, the maps $\rho\circ f_d\circ p_k$ can distinguish elements which are not distinguished by Burau representation.
\end{example}

\section{A new invariant of pure braids: quadruples of points on the same circle}

In this section, we deal with the subgroup ${\widetilde {PB}}_{n}\in PB_{n}$ consisting of those
pure $n$-strand braids for which any two strands have linking number $0$
(in particular, this group contains all Brunnian $n$-strand braids).

Unlike the previous section, we shall not use parameter $d$.

We adopt the above notation: for a braid $\beta=\{\beta_k\}_{k=1}^n$ where $\beta_k\colon[0,1]\to \R^2=\mathbb C$.

Fix two distinct indices $k,l\in \{1,\cdots, n\}$. With a braid $\beta$ we shall associate the cylindrical $(n-2)$-strand braid as follows.

Without loss of generality, we may assume $\beta_{k}(t)\sim 0, \beta_{l}(t)\sim 1 $ (otherwise we just apply fractional linear map taking $\beta_{k}(t)$ and $\beta_{l}(t)$ to $0, 1$ respectively).

If we delete the stands $k,l$ from the braid $\beta$ we get a pure braid in the plane with two punctures
(equivalently, sphere with three punctures $S_{0,3} = {\bar {\mathbb C}}\backslash \{0,1,\infty\}$. Denote
this braid by $q_{k,l}$.

\begin{proposition}
The map $q_{k,l}$ is a homomorphism from the pure braids on $n$ strands ${\widetilde {PB}}_n$ to the group of pure braids on $n-2$ strands $PB_{n-2}(S_{0,3})$.
\end{proposition}

\begin{remark}
The map $q_{k,l}$ can be given by the formula
\begin{equation}
q_{k,l}(\beta)_i(t)=\frac{\beta_i(t)-\beta_k(t)}{\beta_l(t)-\beta_k(t)},\quad i\ne k,l,\ t\in[0,1].
\end{equation}
\end{remark}

To go further, we shall need to formalise the new class of braids, which we call {\em flat-virtual braids}.
By a {\em flat-virtual braid} we mean an equivalence class of braid diagrams having three type of crossings, {\em classical ones},
{\em flat ones}, and {\em virtual ones}, modulo the standard relations: classical braid relations for three classical
crossings, far commutativity relations for any two types of crossings, second Reidemeister moves for two flat crossings,
second Reidemeister move for two virtual crossings, third Reidemeister move with three virtual
crossings and mixed third Reidemeister moves where virtual strand goes through a classical (resp., flat) crossing
or flat strand goes through a classical crossing.

The corresponding groups will be denoted by $FVB_{n}$, the subgroup of pure flat-virtual braids
(kernel of the map to the permutation group) will be denoted by $PFVB_{n}$. By forgetting the under-overcrossing structure at the classical
crossings, we get obvious maps $FVB_{n}\to FB_{n}, PFVB_{n}\to PFB_{n}$ to
the corresponding flat braid groups.

We may think of flat-virtual braids just in terms of crossings which are not virtual: we just
put classical crossings (with information about overpasses and underpasses) and flat crossings
(without such information) and indicate the ways to connect them; the connections
will give rise to strands containing virtual crossings which are defined up to detours. Later
on, we shall enumerate the strands globally, i.e., we refer to $\gamma_{i}(t), t\in [0,1]$ as
$i$-th strand. At each flat or classical crossing, one also has to indicate which strand
goes from the north-east to south-west (called NE-strand), the other stand goes from north-west to
south east (called NW-strand).

Given a pure braid $\gamma\in PB_{n-2}(S_{0,3})$ in $S_{0,3} = {\bar C}\backslash \{0,1,\infty\}$,
we will be interested in those moments $t$ when some two points $\gamma_{i}(t),\gamma_{j}(t)$ belong
to a circle (or line) passing through two punctures.

By a slight perturbation we can make it {\em generic} with respect to such moments (i.e., there are
finitely many such moments $t'$, at each such moment $t'$ there is exactly one pair of indices when
such event happen and this happens transversally).

Given a moment $t_{0}$ where $(0,1,\gamma_{i}(t),\gamma_{j}(t_{0}))$ belong to the same circle or line.
By applying a family of linear-fractional transformations for $t\in [t_{0}-\eps, t_{0}+\eps]$ for $\eps$
small enough, we may assume that $\gamma_{i}(t)\sim \frac{1}{2}$ for these values of $t$. The transversality
here means that $\gamma_{j}(t)$ passes transversally through real line at $t=t_{0}$.

We are defining non-virtual crossings of the braid $\psi(\gamma)$ to be constructed.
Now, with each such moment $t_{0}$ we associate a
classical crossing with $j$-th strand going
over $i$-th strand if $\gamma_j(t_0)\in [\frac{1}{2},1]$,
a classical crossing with $j$-th strand going under $i$-th strand
if $\gamma_j(t_0)\in [0,\frac{1}{2}]$, or a flat crossing if
$\gamma_j(t_0) \in \R \backslash [0,1]$. More specifically, $j$-th strand
is NE at the crossing above if and only if the derivative $\gamma'_j(t_{0})$ has negative
imaginary part.

\begin{figure}[h]
\centering
\includegraphics[width=0.8\textwidth]{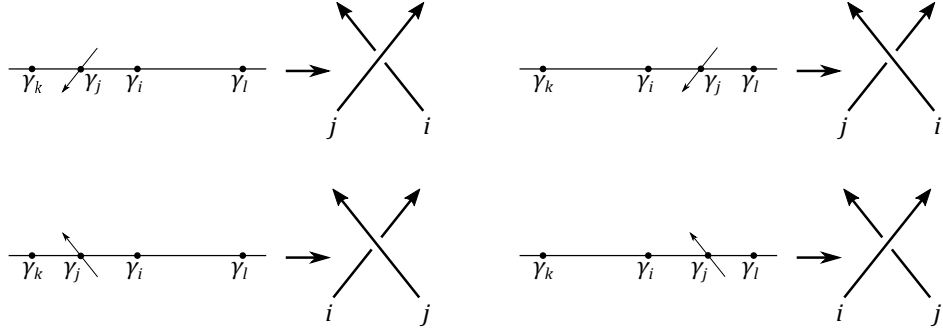}
\caption{Classical crossings of the braid $\psi(\gamma)$}\label{fig:classical_crossings}
\end{figure}

\begin{proposition}
The map $\psi$ defined above is a well defined homomorphism $\psi:PB_{n-2}(S_{0,3})\to PFVB_{n-2}$.
\end{proposition}

\begin{remark}
Choosing another configuration for the strands $k$ and $l$: $\beta_{k}(t)\sim 0, \beta_{l}(t)\sim \infty$ we can rewrite the exceptional moments considered above as follows: we place a classical crossings on the strands $i$ and $j$ if the number $\beta_{i}(t)/\beta_{j}(t)$ is real and positive, and place a flat crossings on the strands $i$ and $j$ if the number $\beta_{i}(t)/\beta_{j}(t)$ is real and negative.

Such description admits the following generalisation. Fix a natural number $d$, and place a classical crossings on the strands $i$ and $j$ if the argument $arg(\beta_{i}(t)/\beta_{j}(t))=0$, and place a flat crossings if $arg(\beta_{i}(t)/\beta_{j}(t))=\frac{2\pi p}d$, $p=1,\dots,d-1$. Thus, we get a series of homomorphisms $\psi_d:PB_{n-2}(S_{0,3})\to PFVB_{n-2}$ where $\psi_2=\psi$.

When $d>2$ we must add the third Reidemeister move on three flat crossings to the relations of $PFVB_{n-2}$.
\end{remark}

We can rewrite the flat-virtual braid $\psi(q_{k,l}(\beta))$ in terms of classical, flat, and virtual generators $\sigma_{i}, \pi_{i},\tau_{i}$, $1\le i\le n-3$.

The relations are
\begin{itemize}
\item $x_iy_j=y_jx_i$ for any $|i-j|>1$ and $x,y\in\{\sigma,\pi,\tau\}$;
\item $x_ix_{i+1}x_i=x_{i+1}x_ix_{i+1}$ for any $i$ and $x\in\{\sigma,\tau\}$;
\item $\pi_i^2=\tau_i^2=1$  for any $i$;
\item $x_ix_jy_i=y_jx_ix_j$  for any $|i-j|=1$, $x\in\{\pi,\tau\}$ and $y\in\{\sigma,\pi\}$.
\end{itemize}

Now, we define a presentation ${\tilde \rho}$ in matrices $M(n-2)$ taking $\sigma_{i},\pi_{i}, \tau_{i}$ to a
block-diagonal matrix with the only non-trivial $2\times 2$-block, where for the corresponding block matrix $\{m_{\cdot, \cdot}\}$
is

\[
\left(
\begin{array}{cc}
   m_{ii}  & m_{i,i+1} \\
   m_{i+1,i} & m_{i+1,i+1}
\end{array}
\right)=
\left(
\begin{array}{cc}
   1-t  & t \\
    1 & 0
\end{array}
\right)
\] for $\sigma_{i}$,

\[
\left(
\begin{array}{cc}
   m_{ii}  & m_{i,i+1} \\
   m_{i+1,i} & m_{i+1,i+1}
\end{array}
\right)=
\left(
\begin{array}{cc}
   0  & s \\
    s^{-1} & 0
\end{array}
\right)
\] for $\pi_{i}$, and

\[
\left(
\begin{array}{cc}
   m_{ii}  & m_{i,i+1} \\
   m_{i+1,i} & m_{i+1,i+1}
\end{array}
\right)=
\left(
\begin{array}{cc}
   0  & r \\
    r^{-1} & 0
\end{array}
\right)
\] for $\tau_{i}$.

By analogy with Proposition~\ref{prop:pHomomorphism} we have the following statement.

\begin{proposition}
The map ${\tilde \rho}$ defines a homomorphism ${\tilde \rho}\colon PFVB_{n-2}\to GL(n-2,\Z[t,t^{-1},s,s^{-1},r,r^{-1}])$.
\end{proposition}

\begin{corollary}
The composition ${\tilde\rho}\circ \psi \circ q_{k,l}\colon {\tilde  {PB}}_n\to GL(n-2,\Z[t,t^{-1},s,s^{-1},r,r^{-1}])$ is a representation of the pure braid group.
\end{corollary}

\subsection*{Aknowledgements}

We are extremely grateful to Huyue Yan for discussing the text and pointing to some remarks. We express our gratitude to Louis H.Kauffman for fruitful and stimulation discussions and Nikita D.Shaposhnik for pointing out some remarks in the first version of the text.


\end{document}